\documentclass[12pt]{article}

\usepackage{amsmath,amsthm}
\usepackage{amssymb}
\usepackage{color}
\usepackage{xspace}
\usepackage[colorlinks=true,
linkcolor=green,
filecolor=brown,
citecolor=green]{hyperref}

\begin{document}
\newtheorem{lemma}{Lemma}
\newtheorem{theorem}{Theorem}
\newtheorem{prop}{Proposition}
\newtheorem{cor}[theorem]{Corollary}
\begin{center}
{\Large
Sierpinski's Triangle and the Prouhet-Thue-Morse Word                           \\ 
}

\vspace{10mm}
DAVID CALLAN  \\
Department of Statistics  \\
\vspace*{-1mm}
University of Wisconsin-Madison  \\
\vspace*{-1mm}
1300 University Ave  \\
\vspace*{-1mm}
Madison, WI \ 53706-1532  \\
{\bf callan@stat.wisc.edu}  \\
\vspace{5mm}

November 18, 2006
\end{center}

\begin{abstract}
   Sierpinski's triangle is a fractal and the Prouhet-Thue-Morse word 
   is sufficiently chaotic to avoid cubes. Here  we observe that there is at 
   least a tenuous connection between 
   them:  the Sierpinski triangle is evident in Pascal's triangle mod 2 
   whose inverse, as an infinite lower-triangular
   matrix, involves the Prouhet-Thue-Morse word.
\end{abstract}

Pascal's triangle mod 2 (Fig. 1b) is a discrete version of the fractal 
known as the Sierpinski triangle \cite{sierpinski}. Left-justified, it 
forms an infinite lower-triangular (0,1)-matrix $S$ with 1s on the 
diagonal (Fig. 1c).

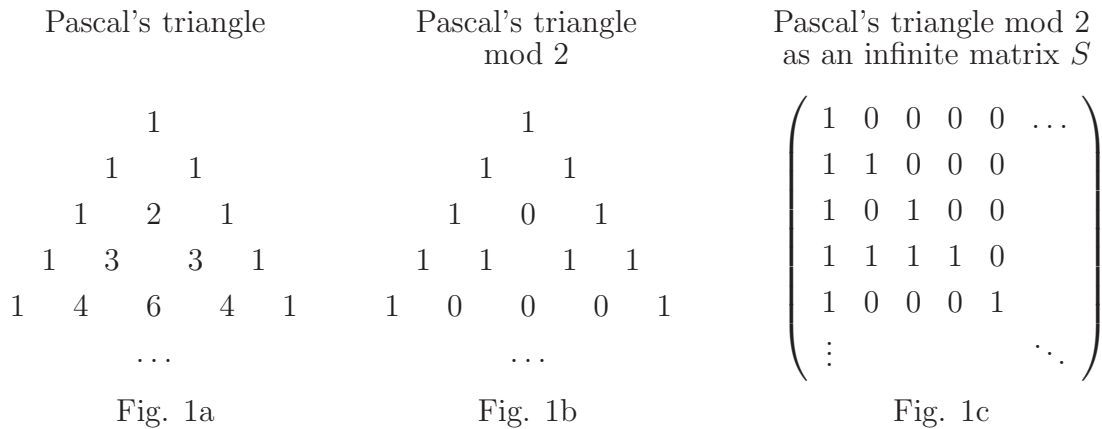
\[
\hspace*{4mm} \stackrel{{ \textrm{\normalsize{Pascal's triangle}}}  }
{\vphantom{mod 2}  } \hspace*{20mm} 
\stackrel{{ \textrm{\normalsize{Pascal's triangle}}}  }
{{ \textrm{mod 2}}  }
\hspace*{18mm}  \stackrel{{ \textrm{\normalsize{Pascal's triangle mod 2}}}  }
{{ \textrm{ as an infinite matrix $S$}}  }
\]
\[
\begin{array}{ccccccccc}
      \!&\!   \!&\!   \!&\!   \!&\! 1 \!&\!   \!&\!   \!&\!   \!&\!    \\
      \!&\!   \!&\!   \!&\! 1 \!&\!   \!&\! 1 \!&\!   \!&\!   \!&\!    \\
      \!&\!   \!&\! 1 \!&\!   \!&\! 2 \!&\!   \!&\! 1 \!&\!   \!&\!   \\
      \!&\! 1 \!&\!   \!&\! 3 \!&\!   \!&\! 3 \!&\!   \!&\! 1 \!&\!    \\
    1 \!&\!   \!&\! 4 \!&\!   \!&\! 6 \!&\!   \!&\! 4 \!&\!   \!&\! 1 \\
      \!&\!   \!&\!   \!&\!   \!&\! \ldots \!&\!   \!&\!   \!&\!   \!&\!  
\end{array}
\hspace*{8mm}
\begin{array}{ccccccccc}
      \!&\!   \!&\!   \!&\!   \!&\! 1 \!&\!   \!&\!   \!&\!   \!&\!    \\
      \!&\!   \!&\!   \!&\! 1 \!&\!   \!&\! 1 \!&\!   \!&\!   \!&\!    \\
      \!&\!   \!&\! 1 \!&\!   \!&\! 0 \!&\!   \!&\! 1 \!&\!   \!&\!   \\
      \!&\! 1 \!&\!   \!&\! 1 \!&\!   \!&\! 1 \!&\!   \!&\! 1 \!&\!    \\
    1 \!&\!   \!&\! 0 \!&\!   \!&\! 0 \!&\!   \!&\! 0 \!&\!   \!&\! 1 \\ 
      \!&\!   \!&\!   \!&\!   \!&\! \ldots \!&\!   \!&\!   \!&\!   \!&\! 
\end{array}
\hspace*{12mm}
\left(\begin{array}{ccccccccc}
    1 & 0 & 0 & 0 & 0 & \ldots \\
    1 & 1 & 0 & 0 & 0 &   \\
    1 & 0 & 1 & 0 & 0 &  \\
    1 & 1 & 1 & 1 & 0 &   \\
    1 & 0 & 0 & 0 & 1 &  \\
    \vdots &   &   &   &   & \ddots
\end{array} \right)
\]
\[
 \textrm{Fig. 1a} \hspace*{35mm} \textrm{Fig. 1b} 
\hspace*{42mm} \textrm{Fig. 1c}
\]

The Prouhet-Thue-Morse word on a 2-letter alphabet $\{a,b\}$ can be 
formed as follows. Start 
with $a$, switch letters and append to get $ab$. Again switch letters 
and append to get $abba$. Repeat to get $abbabaab$ and iterate. The 
result is the infinite Prouhet-Thue-Morse word $t(a,b)$ ($t$ for 
Thue) that crops 
up in diverse contexts \cite{ptm}. Here we 
observe that the word $t$ is related to the matrix $S$.
\begin{theorem}
    $S^{-1}$ is a $(-1,0,1)$-matrix. It has the 
   same pattern of zeroes as $S$ and the nonzero entries in each 
   column form the Prouhet-Thue-Morse word\, $t(1,-1)$. 
\end{theorem}

This is an immediate corollary of a more general result, Theorem 4 below. 

For nonnegative 
integers $i$ and $j$, say $i$ is (binary-)free of $j$ if the binary 
expansion of $i$ has 0s in the positions where $j$ has 1s, 
equivalently, if the operation of adding $i$ and $j$ in base 2 involves 
no ``carries''. By Kummer's well known criterion for the power of a 
prime dividing a binomial coefficient \cite[Ex. 5.36]{gkp}, $i$ is free of 
$j\Leftrightarrow \binom{i+j}{j}$ is odd. Thus the matrix $S$ has 0s precisely in the 
positions $(i,j)$ where $i-j$ is not free of $j$. For nonnegative integer $n$, 
let $b(n)$ denote the sum of the binary digits of $n$. For example, 
$6_{2}=110$ and $b(6)=2$. The Prouhet-Thue-Morse word $t(1,-1)$ has 
the explicit form 
$\big( (-1)^{b(n)} \big)_{n\ge 0}$ \cite{ptm}. With $x$ an 
indeterminate, let $S(x)$ denote the infinite lower-triangular matrix 
defined by
\[
S(x)_{ij}=\begin{cases} 
x^{b(i-j)} & \textrm{if $i\ge j \ge 0$ and $i-j$ is free of $j$, and} \\
0 & \textrm{otherwise.}
\end{cases}
\]
Thus 
\[
S(x)=\left(
\begin{array}{ccccccccc}
    1 &  &  &  &  &  &  &  &   \\
    x & 1 &  &  &  &  &  &  &   \\
    x & 0 & 1 &  &  &  &  &  &   \\
    x^{2} & x & x & 1 &  &  &  &  &   \\
    x & 0 & 0 & 0 & 1 &  &  &  &   \\
    x^{2} & x & 0 & 0 & x & 1 &  &  &   \\
    x^{2} & 0 & x & 0 & x & 0 & 1 &  &   \\
    x^{3} & x^{2} & x^{2} & x & x^{2} & x & x & 1 &   \\
    \vdots &  &  &  &  &  &  &  & \ddots
\end{array}\right)
\]
and $S(1)=S$.
\newpage
\begin{theorem}
   $S(x)S(y)=S(x+y)$.
\end{theorem}
\textbf{Proof}\quad First, as an example, consider $i=47,\,j=9$ 
and $k$ between $j$ and $i$ of the form displayed, where $a,b,c$ are 
bits (0 or 1) and a prime superscript indicates the complementary bit: 
$a'=1-a$.
\begin{center}
\begin{tabular}{cccccccc}
    \textrm{integer} & & \multicolumn{6}{c} {\rule[-2mm]{0mm}{8mm}binary expansion}  \\ 
    \hline 
    $i$ & \phantom{dist} & 1&0&1&1\,&\,1 & 1  \\
    $j$ &  & 0 & 0 & 1 & 0 & 0 & 1  \\
    $i-j$ &  & 1 & 0 & 0 & 1 & 1 & 0  \\
    $k$ &  & $a$ & 0 & 1 & $b$ & $c$ & 1  \\
    $i-k$ &  & $a'$ & 0 & 0 & $b'$ & $c'$ & 0 \\
    $k-j$ &  & $a$ & 0 & 0 & $b$ & $c$ & 0
\end{tabular}
\end{center}
Here, $i-j$ is free of $j$. Also, $k$ has 1s where $j$ has 1s and 0s where 
$i$ has 0s. And so $i-k$ is free of $k$ and $k-j$ 
is free of $j$ and these are the only such $k$. In the sum 
$\sum_{k=j}^{i}S(x)_{ik}S(y)_{kj}$ for the $(i,j)$ entry of 
$S(x)S(y)$, $k$ contributes $x^{a'+b'+c'}y^{a+b+c}$ and the sum over 
all bits $a,b,c$ is $(x+y)^{3}$, the $(i,j)$ entry of $S(x+y)$. This 
works in general, as we now demonstrate. 

Suppose $i\ge j$. The $(i,j)$ entry of $S(x)S(y)$ is 
$\sum_{k=j}^{i}S(x)_{ik}S(y)_{kj}$. For $0 \le k \le i$, if $i-k$ is 
free of $k$ then both must have 0s in the positions where $i$ has 0s. 
The bits of $k$ in the positions where $i$ has 1s are arbitrary, and 
then $i-k$ has the complementary bits in these positions. For example, 
with $i=1\,0\,1\,1\,1\,1$ (in binary), $k$ must have the binary form 
$a\,0\,b\,c\,d\,e$ so that $i-k=a'\,0\,b'\,c'\,d'\,e'$. If, further, $k\ge j$ and $k-j$ 
is free of $j$ then the bits of $k$ are further restricted: they must 
be 1 in each position where $j$ has a 1. In short, if $i-k$ is 
free of $k$ and $k-j$ 
is free of $j$, then $k$ must have 0s where $i$ has 0s and 1s where 
$j$ has 1s and is unrestricted where $i$ has a 1 and $j$ has a 0.
In particular, the existence of $k \in [j,i]$ with $i-k$ 
free of $k$ and $k-j$ free of $j$ implies that $i$ must have a 1 in 
each position where $j$ has a 1 and hence $i-j$ is free of $j$. So, 
if $i-j$ is \emph{not} free of $j$, then 
$\big(S(x)S(y)\big)_{ij}=\sum_{k=j}^{i}0=0=S(x+y)_{ij}$.
On the other hand, if $i-j$ is free of $j$, suppose there are $t\ge 
0$ positions where $i$ has a 1 and $j$ has a 0 (and so $b(i-j)=t$). As 
above, the $k \in [j,i]$ for which  $i-k$ is
free of $k$ and $k-j$ is free of $j$ are unrestricted in these 
positions and both $i-k$ and $k-j$ have 0s in all other positions. 
Hence
\begin{eqnarray*}
   \big(S(x)S(y)\big)_{ij} & = &  \sum_{k=j}^{i}S(x)_{ik}S(y)_{kj} \\
     & = & \sum_{(i_{1},\ldots,i_{t})\in \{0,1\}^{t}}x^{i_{1}+\ldots 
     +i_{t}}y^{i_{1}'+\ldots +i_{t}'}  \\
     & = & \sum_{m=0}^{t}\binom{t}{m}x^{m}y^{t-m}  \\
     & = & (x+y)^{t}  \\
     & = & S(x+y)_{ij}.
\end{eqnarray*}  
Induction yields
\begin{cor}
    For $q$ a positive integer, $S(x)^{q}=S(qx)$.
\end{cor}
\begin{theorem}
    For rational $r$, $S^{r}=S(r)$.
\end{theorem}
\textbf{Proof} \quad For $r=p/q$ with $p,q$ positive integers, this 
follows from
\[
S(p/q)^{q}\underset{ \textrm{Cor. 3} }{=}S(p)\underset{ \textrm{Cor. 3} }{=}S(1)^{p}=S^{p}.
\]
For negative $r$, it now suffices to show that $S^{-1}=S(-1)$ and this 
follows from
\[
S(-1)\,S=S(-1)\,S(1)\underset{ \textrm{Thm. 2} }{=}S(0)=I.
\]

\textbf{Added in Proof.\ } Roland Bacher informs me that he has obtained 
these results more simply by observing that the $2^{k}\times 2^{k}$ 
upper left submatrix of $S(x)$ is the $k$-fold Kronecker product of $
  \left( \begin{smallmatrix}
  1 & 0 \\
    x & 1
\end{smallmatrix}\right)$ \cite{sympascal,rec}. Emmanuel Ferrand 
treats similar material in an interesting recent paper \cite{ferrand06}.


\begin{thebibliography}{99}
    
\bibitem{sierpinski} MathWorld, 
\htmladdnormallink {The Sierpinski 
Sieve}{http://mathworld.wolfram.com/SierpinskiSieve.html} .

\bibitem{ptm} J.-P. Allouche and J. O. Shallit, 
\htmladdnormallink{The Ubiquitous Prouhet-Thue-Morse Sequence}{http://www.cs.uwaterloo.ca/~shallit/Papers/ubiq.ps}
in C. Ding. T. Helleseth and H. Niederreiter, eds., Sequences and Their Applications: 
Proceedings of SETA '98, Springer-Verlag, 1999, pp. 1-16. 

\bibitem{gkp} Ronald L. Graham, Donald E. Knuth, Oren Patashnik, 
\emph{Concrete Mathematics} 
(2nd edition), Addison-Wesley, 1994.

\bibitem{sympascal} Roland Bacher and Robin Chapman, Symmetric Pascal matrices 
modulo p, \emph{European J. Combin.} \textbf{25} (2004), 459--473,
\htmladdnormallink {http://front.math.ucdavis.edu/math.NT/0212144}{http://front.math.ucdavis.edu/math.NT/0212144} .

\bibitem{rec} Roland Bacher, La suite de Thue-Morse et la cat\'{e}gorie 
\textbf{Rec},
\emph{Comptes Rendues Acad. des Sci. Paris} Ser I, \textbf{342} (2006), 
161--164.

\bibitem{ferrand06} Emmanuel Ferrand, An analogue of the Thue-Morse 
sequence, preprint, 2006,
\htmladdnormallink 
{http://www-fourier.ujf-grenoble.fr/\raisebox{-2pt}{\~{}}eferrand/ZTM.pdf}{http://www-fourier.ujf-grenoble.fr/~eferrand/ZTM.pdf} .






\end{thebibliography}
\end{document}